\documentclass{amsart}
\usepackage{graphicx}
\usepackage{psfrag}
\usepackage{amsthm}
\usepackage{amscd}
\usepackage{amsfonts}
\usepackage{amstext}
\usepackage{amssymb}
\usepackage{amsmath}
\usepackage{amsxtra}
\usepackage{plain}
\usepackage[all]{xy}

\newcommand{\comm}[1]{}

\numberwithin{equation}{section}
\let\newpf\proof  
\newenvironment{pf}{\newpf{\bf proof }}{\qed\endtrivlist}

\let\newpf\proof

\def\be{\begin{equation}}
\def\ee{\end{equation}}

\def\ba{{\begin{align}}}
\def\ea{{\end{align}}}

\def\bp{{\begin{prop}}}
\def\ep{{\end{prop}}}

\def\bt{{\begin{thm}}}
\def\et{{\end{thm}}}

\def\bpf{{\begin{pf}}}
\def\epf{{\end{pf}}}

\def\ba*{{\begin{align*}}}
\def\ea*{{\end{align*}}}

\def\bg{{\begin{gather}}}
\def\eg{{\end{gather}}}

\def\bg*{{\begin{gather*}}}
\def\eg*{{\end{gather*}}}

\def\bex{\begin{exercise}}
\def\eex{\end{exercise}}
\newtheorem{prop}{Proposition}
\newtheorem{cor}[prop]{Corollary}

\newtheorem{lem}[prop]{Lemma}
\newtheorem{thm}[prop]{Theorem}

\newtheorem{rem}{\bf Remark}

{}

{}
{}
{}

\newcommand{\Z}{\mathbb{Z} }
\newcommand{\N}{\mathbb{N} }

\newcommand{\R}{\mathbb{R} }

\newcommand{\tr}{\operatorname{tr}}

\newcommand{\End}{\rm{End}}

\newcommand{\Tr}{\operatorname{Tr}}

\newcommand{\Trg}{\operatorname{Tr_{\langle g\rangle}}}
\newcommand{\lgr}{\operatorname{\langle g\rangle}}
\newcommand{\ler}{\operatorname{\langle e\rangle}}

\newcommand{\T}{\operatorname{T}}

\newcommand{\tKk}{\operatorname{\tilde K\sp k}}

\newcommand{\tnab}{\operatorname{\tilde\triangle}}
\newcommand{\tnabk}{\operatorname{\tilde\triangle\sp k}}
\newcommand{\tnabks}{\operatorname{\tilde\triangle_s^k}}

\begin{document}
\title{Delocalized Betti numbers and Morse type inequalities}
\author{M.E. Zadeh}
\address{Mathematisches Institut
Georg-August-Universität
Göttingen, Germany\and Institute for Advanced Studies in Basis Sciences, Zanjan, Iran}
\email{zadeh@uni-math.gwdg.de}
\footnotetext{\emph{Key words and phrases}:Novikov-Shubin inequalities,
Delocalized Betti numbers, Witten's Laplacian.}
\footnotetext{\emph{AMS 2000 Mathematics Subject Classification}:
primary 58J35, secondary 58E05}
\maketitle
\begin{abstract}
In this paper we state and prove Morse type inequalities for Morse
functions as well as for closed differential 1-forms. These
inequalities involve delocalized Betti numbers. As an immediate
consequence, we prove the vanishing of delocalized Betti numbers of
manifolds fibering over the circle.
\end{abstract}

\section {Introduction}
Given a manifold $M$ and a real Morse function $f$ on $M$ the
following Morse inequalities
establish relations between the very basic topology of $M$ and the number
of critical points of order $j$ denoted by $C_j$ (cf. \cite{Milnor-Morse})
\[C_k-C_{k-1}+\dots\pm C_0\geq\beta^k-\beta^{k-1}+\dots\pm \beta^0\]
Here $\beta^j:=\dim H^j(M,\R)$ is the $j$-th Betti number of $M$.
These relations have been
the subject of many significant generalizations. S.Novikov and
M.Shubin have proved in \cite{NovikovShubin-Morse} that
these inequalities hold if the Betti numbers being replaced by the
$L^2$-Betti numbers. The $L^2$-Betti numbers(or von Neumann Betti numbers)
were introduced by M.Atiyah in \cite{Atiyah-von} in his investigation on
equivariant index theorem.
The Morse theory for closed 1-forms has been introduced
by S.Novikov and he has proved in \cite{Novikov-Hamiltonian} that the
Morse inequalities can be generalized to closed 1-forms if one replace the
Betti
numbers by the Novikov numbers.
In \cite[theorem 1]{Farber-Von} it is shown that the Novikov-Shubin inequalities
hold as well for closed 1-forms. In this paper we are interested in the
delocalized Betti numbers which are introduced by J.Lott in \cite{Lott-delocalized}.
These delocalized
Betti numbers are not yet well studied and enjoy properties which are not satisfied by
the ordinary or
$L^2$-Betti numbers, e.g the delocalized Betti numbers of any manifold with
free abelian
fundamental group and of hyperbolic manifolds vanish. However In this paper
we show that
some appropriate combinations of delocalized Betti numbers satisfy the Novikov-Shubin
and the Farber inequalities (see theorems \ref{ahgar} and \ref{delahgar}).
We prove the delocalized Novikov-Shubin inequalities by following
J.Roe's account\cite{Roe-elliptic}
of the E.Witten approach to Morse theory \cite{Witten-supersymmetry}.
As a consequence of this
method we reprove the vanishing of the delocalized Euler character of $M$.
To prove the delocalized Morse inequalities for closed 1-forms
we use these inequalities for Morse function and follow the method used by M.Farber.
As a consequence of the Morse inequalities for closed 1-forms we prove that the
delocalized Betti numbers of manifolds which fiber over the circle
vanish. This vanishing theorem for $L^2$-Betti numbers was conjectured by
M.Gromov and proved by W.L\"uck in \cite{Luck-L2}.

In section \ref{section1} we state and prove the Morse inequalities in a
very general analytic framework. And in section \ref{section2} we use these
inequalities to prove the Morse inequalities for Morse functions and delocalized
Betti numbers. In these two section we follow closely the methods used by J.Roe in
\cite[chapter 14]{Roe-elliptic}. We use the main theorem of section \ref{section2}
to prove the Morse inequalities for closed 1-forms and delocalizd Betti numbers in section \ref{section3}.

\section{General analytic Morse inequalities}\label{section1}

Let $(M,g)$ be a closed oriented riemannian manifold and denote
by $\triangle\sp k$ the laplacian operator acting on differential $k$-forms on $M$.
Let $(\tilde M,\tilde g)$ be the riemannian universal covering of $M$ with $\tilde g=\pi^*g$,
where $\pi$ is the covering map. We denote by $G$ the fundamental group of $M$ and by  $\tnabk=d^* d+d d^*$ the laplacian operator
acting on $L^2$-elements of $\Omega^k(\tilde M)$.
For  $0\leq k\leq n$ let $\T_k$ be a real valued non-negative
and continuous trace on the space of all smoothing $G$-invariant operators
on $L^2(\tilde M,\Lambda^kT\tilde M)$. The continuity is
understood with respect to the uniform convergence of Schwartz kernels on compact
subsets of $M\times M$. We usually omit the subfix $k$ and denote these
traces by the same symbol $\T$.
Let $\tilde P^k$ denotes the orthogonal projection on $\ker\tnabk$ which is a
smoothing operator on $L^2(\tilde M,\Lambda^kT\tilde M)$, c.f. \cite{Atiyah-von}.
We define the $k$-th T-Betti number by the following relation
\begin{equation}
\beta_{\T}^k:=\T(\tilde P^k).
\end{equation}
For $\phi$ a rapidly decreasing non-negative smooth function on $\R^{\geq0}$ with
$\phi(0)=1$. The operator $\phi(\tnabk)$ is a smoothing operator and so we
can define $\mu_{\T}^k=\T(\phi(\tnabk))$.
Notice that $\beta_T^k$ and $\mu_T^k$ are both non-negative real numbers.
\begin{thm}[analytic Morse inequalities]\label{anamorse}
For $0\leq k\leq n=\dim M$ we have the following inequalities
\[\mu^k_{\T}-\mu^{k-1}_{\T}+\dots\pm\mu^0_{\T}\geq
\beta_{\T}^k-\beta_{\T}^{k-1}+\dots\pm \beta_{\T}^0,\]
and the equality holds for $k=n$.
\end{thm}
\begin{pf}
Let $\{\phi_m\}_m$ be a sequence of non-negative rapidly decreasing
smooth functions on $\R^{\geq0}$ which converges to zero outside of
$0\in\R$ and $\phi_m(0)=1$. The operator $\phi_m(\tnabk)$ consists
of non-negative smoothing operators with smooth Schwartz kernel
kernel $K_m$. The sequence of kernels $K_m$ converges uniformly on
compact subsets of $M\times M$ to the kernel $\tilde K$ of the
projection $\tilde P^k$. In fact the spectral theorem for
self-adjoint operator $\tnabk$ (see, e.g. \cite[theorem VIII.
5]{Reed&Simon-1}) implies that \(\phi_m(\tnabk)\omega\to\tilde P^k\)
for any $L^2$-differential $k$-form $\omega$, i.e.
$\phi_m(\tnabk)\to\tilde P^k$ weakly in $\mathcal L(\mathcal
D,\mathcal D')$. Since the strong and the weak topology of $\mathcal
L(\mathcal D,\mathcal D')$ coincide on bounded subsets, e.g.
converging sequence one conclude the convergence
$\phi_m(\tnabk)\to\tilde P^k$ in the topological space $\mathcal
L(\mathcal D,\mathcal D')$ with strong topology. The Schwartz kernel
theorem asserts that this topological space is isomorphic to
$\mathcal D'(\tilde M\times\tilde M, \Lambda^kT^*\tilde M\times
\Lambda^kT^*\tilde M)$. and with respect to this isomorphism the The
convergence $\\phi_m(\tnabk)\to\tilde P^k$ reads the convergence of
the kernels $K_m\to\tilde K$ which implies the assertion. Therefore
by continuity of $T$ and by finiteness of $|\lgr|$ we obtain
\begin{equation}\label{nakar}
\beta_T^k=\lim_{m\to\infty}\T(\phi_m(\tnabk)).
\end{equation}

The function $\phi-\phi_m$ is non-negative, rapidly decreasing and vanishing
at $0$, so it takes the form $x\psi_m^2(x)$ where $\psi_m$ is non-negative and
rapidly decreasing. We have
\begin{align*}
\T(dd^*\psi_m(\tnabk)^2)&=\T(\psi_m(\tnabk)d\,d^*\psi_m(\tnabk))\\
&=\T(d^*\psi_m(\tnabk)^2\,d)\\
&=\T(d^*d\psi_m(\tnab^{k-1})^2).
\end{align*}
In the last step we have used the commutation relation
$\tnabk\,d=d\,\tnab^{k-1}$. Therefore
\[(\mu_{\T}^k-T(\phi_m(\tnabk))-(\mu_{\T}^{k-1}-T(\phi_m(\tnab^{k-1}))+
\dots\pm (\mu_{\T}^0-T(\phi_m(\tnab^0))=T(d^*d\,\psi_m(\tnabk)^2))~.\]
Since $d^*d\,\psi_m(\tnabk)^2$ is a non-negative smoothing operator,
the right side of the above equality is non-negative for $k<n$ and is zero for $k=n$.
By tending $m$ toward infinity and using relation \eqref{nakar} we get the desired
inequalities.
\end{pf}

\section{Delocalized Novikov-Shubin inequalities}\label{section2}
We recall the definition of \emph{delocalized Betti numbers} as they are introduced
by J. Lott in \cite{Lott-delocalized}. We keep the notation of the previous section.
Let $\tilde P$ be a smoothing $G$-invariant operator acting on $\Lambda^kT^*\tilde M$
with Schwartz  kernel $\tilde K$ which is rapidly decreasing far from the diagonal of $\tilde M\times\tilde M$. For each $\tilde x\in\tilde M$ and $h\in G$
one can identify both $\Lambda^kT_{\tilde x}^*\tilde M$ and
$\Lambda^kT_{h.\tilde x}^*\tilde M$ with $\Lambda^kT_x^* M$ where $x=\pi(\tilde x)$.
Hencefore $\tilde K(\tilde x,h.\tilde x)$ can be considered as an element of
$\End(\Lambda^kT_x^*M)$.  Consequently, given a conjugacy class $\lgr$ of $G$ the sum
$\sum_{h\in\lgr}\tKk(h.\tilde x,\tilde x)$ is finite and,
as a function of $\tilde x\in\tilde M$, is invariant with respect to
the action of $G$. Therefore this function
pushes down and defines a smooth section $K_{\lgr}$ of the bundle
$\End(\Lambda^k(T^*M))$ over $M$.
The following relation defines the delocalized trace $\Tr_{\lgr}$
\begin{equation}\label{deltrace}
\Trg(\tilde P):=\int_M\tr K_{\lgr}\,d\mu_g.
\end{equation}
As in the previous section let $\tnabk$ denotes the laplacian operator
acting on $L^2$-sections of $\Lambda^kT\tilde M$.
The orthogonal projection $\tilde P^k$ on $\ker(\tnabk)$ is a smoothing
$G$-invariant operator on
$L^2(\tilde M,\Lambda^kT\tilde M)$ with kernel $\tilde K^k$ (see, e.g. \cite{Atiyah-von}).
The $k$-th delocalized Betti number $\beta_{\lgr}^k$ is defined as follows
\begin{equation}
\beta_{\lgr}^k=\Tr_{\lgr}(\tilde P^k).
\end{equation}
Equivalently if $\phi_m$ is a sequence of real functions as in the theorem
\ref{anamorse}, e.g. $\phi_m(x)=e^{-x/m}$ then  by \eqref{nakar}
\begin{equation}\label{nakdar}
\beta_T^k=\lim_{m\to\infty}\Tr_{\lgr}(\phi_m(\tnabk)).
\end{equation}
\begin{rem}
Notice that the trace $\Tr_{\lgr}$ is not always finite since it is not known
weather the kernel $\tilde K_{\lgr}$ is rapidly decreasing far from the diagonal
of $\tilde M\times \tilde M$.
Nevertheless if $\lgr$ is a finite set then the above delocalized Betti
numbers are well defined. This is why we restrict ourself, from now on,
to the class $\mathcal C(G)$ of finite conjugacy classes.
\end{rem}

Since the delocalized traces $\Tr_{\lgr}$ are not in general positive,
we cannot apply the theorem \ref{anamorse} to these traces
(except for $k=n$). Instead
we consider the following  linear functional which are introduced in
\cite{Knebusch-approximation} and are proved to be positive traces
\begin{gather}
\T_{\lgr}:=\Tr_{\ler}+\frac{1}{|\lgr|}\Tr_{\lgr}\label{trreal}
\end{gather}
\begin{lem}
The linear functional $\T_{\lgr}$ is a positive trace on the space of
$G$-invariant smoothing operator.
\end{lem}
\begin{pf}
For each $g\in G$ the linear functional $\Tr_{\lgr}$ is a trace
(cf, \cite[lemma 2]{Lott-delocalized}) so $\T_{\lgr}$ is also a trace.
In below we show that it is positive, i.e. $\T_{\lgr}(\tilde P)\geq0$ if
$\tilde P=\tilde Q^*\tilde Q$ where $\tilde Q$ is a $G$ invariant
smoothing operator on $H:=L^2(\tilde M,\Lambda^kT^*\tilde M)$.
Let $\{\theta_k\}_{k\in \N}$ be an orthonormal basis for
$H$ and for $h\in G$ put
$(h.\theta)(\tilde x):=\theta(h.\tilde x)$. We identify $M$ with a
distinguished fundamental domain of $G$ in $\tilde M$. With this identification we have

\begin{align*}
\Tr_{\ler}(\tilde P)
&=\sum_k\int_M \langle\tilde Q\theta_k\,(\tilde x),
\tilde Q\theta_k\,(\tilde x)\rangle\,d\mu_{\tilde g}\\
&=(\sum_k\int_M \langle\tilde Q\theta_k\,(\tilde x),
\tilde Q\theta_k\,(\tilde x)\rangle\,d\mu_{\tilde g})^{\frac{1}{2}}.
(\sum_k\int_M \langle\tilde Q\theta_k\,(\tilde x),
\tilde Q\theta_k\,(\tilde x)\rangle\,d\mu_{\tilde g})^{\frac{1}{2}}\\
&=(\sum_k\int_{M} \langle\tilde Q\theta_k\,(\tilde x),
\tilde Q\theta_k\,(\tilde x)\rangle\,d\mu_{\tilde g})^{\frac{1}{2}}.
(\sum_k\int_{M} \langle\tilde Q\theta_k\,(h.\tilde x),
\tilde Q\theta_k\,(h.\tilde x)\rangle\,d\mu_{\tilde g})^{\frac{1}{2}}\\
&\geq \sum_k(\int_{M} \langle\tilde Q\theta_k\,(\tilde x),
\tilde Q\theta_k\,(\tilde x)\rangle\,d\mu_{\tilde g})^{\frac{1}{2}}
(\int_{M} \langle\tilde Q\theta_k\,(h.\tilde x),
\tilde Q\theta_k\,(h.\tilde x)\rangle\,d\mu_{\tilde g})^{\frac{1}{2}}
\end{align*}
Here the third equality follows from the $G$-invariance of the
riemannian volume element, while the last inequality is the
Cauchy-Schwarz inequality for $L^2$-sequences. The following pairing
is a symmetric bilinear function on $H$
\[\langle\omega,\eta\rangle:=\int_M\langle\omega(\tilde x),
\eta(\tilde x)\rangle\,d\mu_{\tilde g}\] By applying the
Cauchy-Schwarz inequality arising from this bilinear function to the
last expression in above, we get the following inequality
\[\Tr_{\ler}(\tilde P)\geq |\sum_k\int_{M} \langle\tilde Q\theta_k\,(\tilde x),
\tilde Q\theta_k\,(\tilde h.x)\rangle\,d\mu_{\tilde g}|.\]
This inequality beside the following relation
\begin{equation*}
\Tr_{\lgr}(\tilde P)=\sum_{h\in\lgr}\sum_k\int_{M}
\langle\tilde Q\theta_k\,(\tilde x),
\tilde Q\theta_k\,(\tilde h.x)\rangle\,d\mu_{\tilde g}
\end{equation*}
implies the following relation which prove the assertion of the lemma
\[|\lgr|\Tr_{\ler}(\tilde P)\geq|\Tr_{\lgr}(\tilde P)|~.\]
\end{pf}
Using the trace $\T_{\lgr}$ we define the following combination of delocalized Betti numbers
\begin{gather*}
\gamma_{\lgr}^k:=\T_{\lgr}(\tilde P^k)=\beta_{\ler}^k+\frac{1}{|\lgr|}\beta_{\lgr}^k~.
\end{gather*}

Now let $f$ be a Morse function on $M$ and denote by $\tilde f$ its lifting to $\tilde M$.
For $s>0$ put $d_s:=e^{-s\tilde f}de^{s\tilde f}$ and $d^*_s:=e^{s\tilde f}d^*e^{-s\tilde f}$
and define the Witten deformed laplacian, acting on $L^2$-elements of $\Omega^k(\tilde M)$, by the following relation
\[\tnab_s^k:d_sd_s^*+d_s^*d_s\]
This deformed laplacian is a perturbation of the laplacian $\tnabk$ by
differential operators of order zero. So it is an elliptic second order differential
operator and its $L^2$-kernel consists of smooth differential
$k$-forms and the projection on this kernel is a smoothing operator, c.f.
\cite{Atiyah-von}. Just as in above
we can define the deformed  $k$-th delocalized Betti number $\beta_{\lgr}^k(s)$.
In fact these deformed delocalized Betti numbers $\beta_{\lgr}^k(s)$ are
independent of $s$.
To see that consider the deformed de-Rham complex of $L^2$-differential forms
\[\cdots\longrightarrow\Omega^k(T^*\tilde M)\stackrel{d_s^{k}}{\longrightarrow}
\Omega^{k+1}(T^*\tilde M)\longrightarrow\cdots\] This complex is
equivariant with respect to the action of the group $G$, so the
$k$-th homology vector space of this complex $H_s^k(\tilde M,\R)$ is
a $G$-vector space. This real $G$-vector space is isomorphic to the
kernel of  $\tnabks$. The isomorphism associates to an element in
$\ker\tnabks$ its class in $H_s^k(\tilde M,\R)$ and is clearly
$G$-equivariant. Moreover the above deformed complex is isomorphic
to the ordinary de-Rham complex (corresponding to $s=0$) through
conjugation with $e^{s\tilde f}$ which is $G$-equivariant as well.
We conclude that there is a $G$-equivariant isomorphism between
$\ker\tnabks$ and $\ker\tnabk$. Since the action of $G$ on
$L^2$-differential form is symmetric this implies that the
orthogonal projections on $\ker\tnabks$ and on $\ker\tnabk$ are
$G$-similar. Therefore by taking the $\Tr_{\lgr}$ we conclude that
the deformed Betti numbers are actually independent of $s$. In
particular
\begin{equation}\label{invbet}
\gamma_{\lgr}^k(s)=\gamma_{\lgr}^k~;\hspace*{1cm}0\leq k\leq n.
\end{equation}
These relations and the theorem \ref{anamorse}  applied to the deformed laplacian
implay the following inequalities (equality holds for $k=n$)
\begin{equation}\label{moine}
\mu^k_{\lgr}(s)-\mu^{k-1}_{\lgr}(s)+\dots\pm\mu^0_{\lgr}(s)\geq
\gamma_{\lgr}^k-\gamma_{\lgr}^{k-1}+\dots\pm \gamma_{\lgr}^0,
\end{equation}
where $\mu_{\lgr}^k(s):=\T_{\lgr}\phi(\tnab_s^k)$. We shall to study
the behavior of $\mu^k_{\lgr}(s)$ when $s$ goes to infinity in order
to prove the following theorem.

\begin{thm}[delocalized Novikov-Shubin inequalities]\label{ahgar}
Let $f$ be a Morse function on $M$ and denote by $C_k$ the number of critical
point of Morse index $k$.
For $0\leq k\leq n=\dim M$ we have the following inequalities
\begin{equation}\label{morc}
C_k-C_{k-1}+\dots\pm C_0\geq\gamma_{\lgr}^k-\gamma_{\lgr}^{k-1}+
\dots\pm \gamma_{\lgr}^0,
\end{equation}
and the equality holds for $k=n$.
\end{thm}

\begin{pf}
At first we recall from \cite{Roe-elliptic} that the deformed laplacian
has the following form
\begin{equation}\label{moz}
\tnab_s^k=\tnabk+sL_0+s^2|d\tilde f|^2
\end{equation}
where $L_0$ is a zeroth order operator and $|d\tilde f|$ is the endomorphism given by the
multiplication by $|d\tilde f|$.
There is a positive constant $C$ such that $|d\tilde f(\tilde x)|\geq C$ when $\tilde x$
is outside the union $\tilde U_r$ of the $r$-neighborhoods of
critical points of $\tilde f$. Here $r$ is sufficiently small so that the
$4r$-neighborhoods of critical points are disjoint.
For the proof of the following lemma we refere to \cite[lemma 12.10]{Roe-elliptic}
where the proof is given for compact manifolds but remain true for non-compact
manifolds as well.

\begin{lem}\label{outdia}
Let $\phi$ be a rapidly decreasing
function as above such that the fourier transform of the function $\psi$ defined
by $\psi(t):=\phi(t^2)$ is supported in $(-r,r)$.
Let $\tilde K$ denote the kernel of the smoothing operator $\phi(\tnab_s^k)$.
Then $\tilde K(.,.)$ tends uniformly to zero on
$\tilde M\verb+\+\tilde U_{2r}\times\tilde M\verb+\+\tilde U_{2r}$ when $s$
goes to infinity.
\end{lem}
On the other hand the Schwartz kernel of $\phi(\tnabks)$ is
supported in the distance $r$ of the diagonal of $\tilde
M\times\tilde M$. So if $\omega$ is a differential $k$-form which is
supported within the distance $2r$ of a critical point of $\tilde f$
then $\phi(\tnabks)\omega$ is supported within the distance $3r$ of
the same critical point. To see this let $\tilde D_s:=d_s+d_s^*$,
then $\tilde\nabla_s=(\tilde D_s)^2$ and by the condition on the
support of $\hat\psi$ we have
\begin{gather}\label{fghtu}
\phi(\tnabks)\,\omega(\tilde x)=\psi(\tilde D_s)\,\omega(\tilde x)
=\int_{-r}^r\hat{\psi}(t)e^{-it\tilde D_s}\omega\,(\tilde x)\,dt
\end{gather}
This relation and the unite propagation speed property for the
deformed Dirac operator $\tilde D_s$ imply that
$\phi(\tnabks)\omega$ is supported within distance $3r$ of the
support of $\omega$. Since the action of $G$ on $\tilde M$ is
uniformly properly discontinuous, for $r$ sufficiently small and for
a non-trivial element $h\in G$, the element $(\tilde x,h.\tilde x)$
is not in the distance $2r$ of the diagonal of $\tilde M\times\tilde
M$. Consequently the previous lemma and the above discussion show
that for a non-trivial conjugacy class $\lgr$ the delocalized traces
$\Tr_{\lgr}$ has no contribution in the value of $\lim
\T_{\lgr}\phi(\tnabks)$ when $s$ goes to infinity, so
\begin{equation}\label{fpo}
\lim_{s\to\infty}\mu_{\lgr}^k(s)=\lim_{s\to\infty}\Tr_{\ler}\phi(\tnabks).
\end{equation}
Now we shall to prove the following equality which prove the desired inequalities of
the theorem
\begin{equation}\label{talop}
 \lim_{s\to\infty}\Tr_{\ler}\phi(\tnabks)=C_k~.
\end{equation}
Here $C_k$ is the number of the critical points of $f$ with index $k$.
For this purpose
Let $\tilde\beta=p^*\beta$ be a $G$-invariant smooth function on $\tilde M$ which is
supported in $\tilde U_{3r}$
and is equal to $1$ on $\tilde U_{2r}$. The above lemma shows that
\begin{equation}\label{loccon}
\lim_{s\to\infty}\Tr_{\ler}\phi(\tnabks)
=\lim_{s\to\infty}\Tr_{\ler}(\tilde\beta\phi(\tnabks))~,
\end{equation}
where $\tilde\beta$ is the pointwise multiplication by
$\tilde\beta$. So, the next step is to study the asymptotic behavior
of $\Tr_{\ler}(\tilde\beta\phi(\tnabks))$ when $s$ goes to infinity.
Since the kernel of $\phi(\tnabks)$ is supported in the distance $r$
of the diagonal, the differential forms which are supported outside
$\tilde U_{4r}$ has no contribution in the value of the expression
in the right hand side of \eqref{loccon}. So to evaluate the value
of this expression we can consider only those differential forms
which are supported in $\tilde U_{4r}$. In fact the previous lemma
and the condition on the support and values of $\beta$ show that we
may consider those differential forms which are supported in $\tilde
U_{3r}$. As for $\tnabks$, the kernel of $\phi(\triangle_s^k)$ is
supported in the distance $r$ of the diagonal of $M\times M$. So for
$r$ sufficiently small one can lift the smoothing operator
$\phi(\triangle_s^k)$ to $\tilde M$.  We have the following equality
\begin{equation}\label{estan}
\phi(\tnabks)\,\omega=p^*\phi(\triangle_s^k)\,\omega~;\hspace{7mm}
\text{supp}(\omega)\subset\tilde U_{3r}~,
\end{equation}
To prove this equality, by relations \eqref{fghtu}, it suffices to show that
\[e^{-it\tilde D_s}\,\omega=p^*(e^{-it D_s})\,\omega~;\hspace{1cm}-r\leq t\leq r~.\]
It is clear that we may assume that the support of $\omega$ is included in a small
ball of radius $3r$. For $|t|\leq r$ the both sides of the above equality define smooth
differential forms
which are solutions of the wave equation \(\partial_t+i\tilde D_s=0\). Moreover
they have the same initial condition $\omega$ for $t=0$. Therefore by uniqueness of wave
operator we get the desired equality. Consequently
\begin{align}
\Tr_{\ler}(\tilde\beta\phi(\tnabks))&=\Tr_{\ler}p^*(\beta\phi(\triangle_s^k))\notag\\
&=\Tr (\beta\phi(\triangle_s^k))\label{bahar}~.
\end{align}
Now the argument leading to relation \eqref{loccon} can be applied to $\triangle_s^k$ to deduce
\[\lim_{s\to\infty}\Tr (\beta\phi(\triangle_s^k))=\lim_{s\to\infty}\Tr \phi(\triangle_s^k)=C_k~.\]
The last equality is the main step of the analytical proof of the Morse inequalities for
the Morse function $f$ on $M$, c.f. \cite[page 192]{Roe-elliptic}. This equality and
the above discussion prove the relation \eqref{talop}. This complete the proof of the theorem.
\end{pf}
The inequality \eqref{morc} can be written in the following form
\begin{equation}\label{moralo}
C_k-C_{k-1}+\dots\pm C_0\geq
\beta_{\ler}^k-\beta_{\ler}^{k-1}+\dots\pm \beta_{\ler}^0+B_{\lgr}^k
\end{equation}
where
\begin{align*}
B_{\lgr}^k&:=\frac{1}{|\lgr|}(\beta_{\lgr}^k
-\beta_{\lgr}^{k-1}+\dots\pm\beta_{\lgr}^0)~.
\end{align*}
Inequalities \eqref{moralo} without the term $B_{\lgr}^k$
at the right hand side are the Morse inequalities for $L^2$-Betti numbers established by S.Novikov and M.Shubin \cite{NovikovShubin-Morse}. The proof of the theorem \ref{ahgar} can be applied to the localized trace $\Tr_{\ler}$ and provides an analytical proof for the Novikov-Shubin inequalities.

\begin{rem}
It is clear that the $k=n$-case equality of the theorem \ref{anamorse} holds even
if the trace $T$ is not real. If we apply this equality to $\Tr_{\lgr}$ and follow
the proof of the theorem \ref{ahgar} we get the vanishing of the delocalized Euler
character
\[\chi_{\lgr}(M):=\beta_{\lgr}^n-\beta_{\lgr}^{n-1}\dots\pm \beta_{\lgr}^0=0.\]
Of course this result can be proved by the heat equation approach to compute the
Euler characteristic.
\end{rem}

\section{Delocalized Morse inequalities for closed 1-forms}\label{section3}
In this section we use the main theorem of the previous section to
prove a delocalized version of the Morse inequalities for closed
1-forms. Following \cite{Farber-Von} we use then these inequalities
to prove the vanishing of the delocalized Betti numbers for sapaces
which fiberat  over the circle. Vanishing of the $L^2$-Betti numbers
of such spaces was conjectured by M.Gromov and proved by W.L\"uck
\cite{Luck-L2}.

Let $\omega$ be a closed differential 1-form on $M$. In a small open
subset $U$ one has
$\omega=df_U$, where $f_U$ is a smooth function on $U$ uniquely
determined up to an additive constant.
A point $p\in U$ is a non-degenerate critical point of $\omega$ with index $j$ if
it is a non-degenerate critical point of $f_U$ with index $j$.
As in the previous section we
denote by $C_j$ the number of these  points. Since $\omega$ is closed the map
\begin{equation}\label{hommap}
\gamma\stackrel{\xi}{\rightarrow}\int_\gamma\omega
\end{equation}
defines a homeomorphism between the fundamental groups $G$ and the additive group
$(\R,+)$.

As in the previous section let $C_j$ denotes the number of critical points
of index $j$.
The following theorem reduces to the theorem \ref{ahgar} if $\omega$
is an exact form.
\begin{thm}[delocalized Morse inequalities for closed 1-forms]\label{delahgar}
For $0\leq k\leq n=\dim M$ the following inequalities hold
\begin{equation}
C_k-C_{k-1}+\dots\pm C_0\geq\gamma_{\lgr}^k-\gamma_{\lgr}^{k-1}+
\dots\pm \gamma_{\lgr}^0,
\end{equation}
and the equality takes place for $k=n$.
\end{thm}
\begin{pf}
For $m\in \Z^{\geq0}$ define the normal subgroup $G_m$ of $G$ by
$G_m:=\xi^{-1}(m\Z)$.
Let $M_m$ denote the corresponding cyclic $m$-sheeted normal covering space of $M$
with $\pi_1(M_m)=G_m$. Denote by $\omega_m$ the lifting of $\omega$ to $M_m$.
We have clearly
\[C_j(\omega_m)=m.C_j(\omega).\]
The relevance of these covering spaces is that the number of critical points
of $\omega_m$ can be approximated by the number of critical points of an exact
form on $M_m$. More precisely there is a constant $C$, independent of $m$,
and an exact Morse 1-form
$\tilde \omega_m$ on $M_m$ such that
\begin{equation}\label{apcrit}
C_j(\omega_m)\leq C_j(\tilde\omega_m)\leq C_j(\omega_m)+C~;~\text{ for } 0\leq j\leq n.
\end{equation}
For the proof of this claim we refere to \cite[section 2]{Farber-Von}.
The Morse theory for
the exact 1-form $\tilde\omega_m=d\tilde f_m$ is exactly the Morse theory
for the Morse function $\tilde f_m$.
If in relation \eqref{deltrace} we were to integrate over the
$m$-sheeted covering $M_m$
of $M$ instead of $M$ then we would get a positive trace say $\T_{\lgr,m}$ which
is actually $m$-times of the trace $\T_{\lgr}$.
The arguments leading to the theorem \ref{ahgar}
can be applied to this situation as well. Since the delocalized Betti
numbers defined by $\T_{\lgr,m}$ are
$m$-times the delocalized Betti numbers defined by $\T_{\lgr}$
the resulting Morse inequalities will be the following ones
with the equality for $k=n$
\[C_k(\tilde\omega_m)-C_{k-1}(\tilde\omega_m)+\dots\pm C_0(\tilde\omega_m)
\geq m.\gamma_{\lgr}^k-m.\gamma_{\lgr}^{k-1}+\dots\pm m.\gamma_{\lgr}^0\]
Using these inequalities and relations \eqref{apcrit} we get
\begin{align*}
\sum_{j=0}^k(-1)^{k-j}C_j(\omega)&=\frac{1}{m}\sum_{j=0}^k (-1)^{k-j}C_j(\omega_m)\\
&\geq\frac{1}{m}\sum_{j=0}^k (-1)^{k-j}C_j(\tilde\omega_m)-\frac{k+1}{m}C\\
&\geq \sum_{j=0}^k (-1)^{k-j}\gamma_{\lgr}^j-\frac{k+1}{m}C
\end{align*}
Now taking the limit when $m$ goes to infinity, we obtain the claime of the theorem.
\end{pf}
As a corollary of the above theorem we show that the delocalized Betti
numbers of a manifold, fibering over the circle, vanish.

\begin{cor}[Vanishing theorem]\label{vanisde}
All delocalized Betti numbers of a manifold fiberating over $S^1$ vanish.
\end{cor}
\begin{pf}
Let $M\stackrel{p}{\rightarrow}S^1$ be a fibration. The pull-back 1-form
$\omega=p^*(d\theta)$ on $M$ has no critical point. Applying the inequalities of
the theorem \ref{delahgar} we obtain the following vanishing result
for $0\leq j\leq n$
\[\gamma_{\lgr}^j(M)=\beta_{\ler}^j+\frac{1}{|\lgr|}\beta_{\lgr}^j=0\]
Since $\beta_{\ler}^j=0$ by the previously mentioned result of W.L\"uck,
we conclude that $\beta_{\lgr}^j=0$.
\end{pf}

\end{document}